\documentclass{scrartcl}

\newif\iffinal
\finalfalse
\finaltrue

\usepackage[showeqnr,final,theoremdefs]{latexdev}

\usepackage[utf8]{inputenc}

\definecolor{myredc}{rgb}{.8,0,0}
\definecolor{mybluec}{rgb}{0,0,.8}
\colorlet{myred}{myredc!70!black}
\colorlet{myblue}{mybluec!50!black}
\colorlet{myorange}{orange!50!white}
\hypersetup{
	pdfauthor={Robert McLachlan, Klas Modin, Olivier Verdier},
}

\usepackage{relsize} 

\usepackage[numbers]{natbib} 

\makeatletter 
\def\@seccntformat#1{\csname the#1\endcsname.\quad} 
\makeatother

\usepackage{titlesec}

\titleformat{\section}
{\color{myred}\normalfont\sffamily\large\bfseries}
{\color{myred}\thesection.}{1em}{}
\titleformat{\subsection}
{\color{myred}\normalfont\sffamily\bfseries}
{\color{myred}\thesubsection.}{1em}{}
\titleformat{\subsubsection}
{\color{myred}\normalfont\sffamily\bfseries}
{\color{myred}\thesubsubsection.}{1em}{}

\titleformat*{\paragraph}{\color{myblue}\bfseries}

\usetikzlibrary{calc}
\usetikzlibrary{matrix}
\tikzstyle{commdiag}=[matrix of math nodes, row sep=3em, column sep=5.5em, text height=1.5ex, text depth=0.25ex,ampersand replacement=\&]
\tikzset{>=stealth}

\usepackage{standalone}

\usepackage{pgfplots}
\pgfplotsset{compat=1.5}

\usepackage{subfig}

\usepackage[affil-sl]{authblk} 

\title{Collective Lie--Poisson integrators on $\R^{3}$}
\author[1]{Robert McLachlan\thanks{\href{mailto:r.mclachlan@massey.ac.nz}{r.mclachlan@massey.ac.nz}}}
\author[2]{Klas Modin\thanks{\href{mailto:klas.modin@chalmers.se}{klas.modin@chalmers.se}}}
\author[3]{Olivier Verdier\thanks{\href{mailto:olivier.verdier@math.uib.no}{olivier.verdier@math.uib.no}}}

\affil[1]{
	Institute of Fundamental Sciences, Massey University, New Zealand
}
\affil[2]{
	Department of Mathematical Sciences, Chalmers University of Technology, Sweden
}
\affil[3]{
	Department of Mathematics, University of Bergen, Norway
}

\date{\today}

\usepackage{amsmath,amssymb}
\usepackage{relsize}

\newcommand{\ee}{\mathrm{e}}
\newcommand{\ii}{\mathrm{i}}
\newcommand{\pair}[2]{\left\langle #1, #2 \right\rangle}
\newcommand{\poisson}[2]{\{ #1, #2 \}}

\providecommand{\norm}[1]{\lVert#1\rVert}
\providecommand{\abs}[1]{\lvert#1\rvert}
\providecommand{\vect}[1]{\boldsymbol{#1}}

\newcommand{\ud}{\mathrm{d}}
\newcommand{\udd}{\,\mathrm{d}}

\newcommand{\R}{{\mathbb R}}
\newcommand{\CC}{{\mathbb C}}

\newcommand{\Diff}{\mathrm{Diff}}
\newcommand{\Xcal}{\mathfrak{X}}

\newcommand*\incl{\imath}

\DeclareMathOperator{\Ad}{Ad}

\newcommand*\SO{\mathrm{SO}}

\newcommand*\SU{\mathrm{SU}}
\newcommand*\su{\mathfrak{su}}

\newcommand{\Rattle}{{\smaller RATTLE}}

\newcommand*\dd{\operatorname{d}}
\newcommand*\SUtwo{\SU(2)}
\newcommand*\sutwo{\su(2)}
\newcommand*\Uone{\mathrm{U}(1)}
\newcommand*\uone{\mathfrak{u}(1)}

\newcommand*\hopf{\pi}
\newcommand*\pback[2]{#1^* #2}
\newcommand*\Uoneact{\Psi_{\theta}}
\newcommand*\Meth[1][h]{\Phi_{#1}}
\newcommand*\meth[1][h]{\varphi_{#1}}
\newcommand*\modifham{F_h}

\begin{document}

\maketitle

\begin{abstract}
	We develop Lie--Poisson integrators for general Hamiltonian systems on~$\R^{3}$ equipped with the rigid body bracket. 
	The method uses symplectic realisation of~$\R^{3}$ on $T^{*}\R^{2}$ and application of symplectic Runge--Kutta schemes.
	As a side product, we obtain simple symplectic integrators for general Hamiltonian systems on the sphere $S^{2}$.


	\textbf{MSC~2010:} 37M15, 65P10, 53D17, 53D20 
\end{abstract}

\section{Introduction} 
\label{sec:introduction}

The general problem of symplectic integration is: given a symplectic manifold $(M,\omega)$ and a Hamiltonian function $H\in C^{\infty}(M)$, construct a symplectic map that approximates the flow $\exp(h X_H)$ of the Hamiltonian vector field $X_H$ on $M$.
All known symplectic integrators require additional structures;
either of the symplectic manifold $(M,\omega)$ or of the Hamiltonian $H$.
The main examples are:
\begin{enumerate}
	\item If $(M,\omega)$ is a symplectic vector space, then symplectic Runge-Kutta methods can be used for arbitrary Hamiltonians (see~\cite[\S\!~VI.4]{HaLuWa2006} and references therein).

	\item If $M=T^*Q$ with the canonical symplectic form, and $Q$ is a submanifold of $\R^{n}$ characterised by $Q= \setc{\vect{q}\in \R^{n}}{\vect{c}(\vect{q})=0}$, then the \Rattle~method can be used for (almost) arbitrary Hamiltonians.
	More generally, \Rattle\ gives symplectic integrators on symplectic manifolds realised as transverse submanifolds defined by coisotropic constraints (see \cite{McMoVeWi2013} and references therein).

	\item If the Hamiltonian function is a sum of explicitly integrable subsystems, then splitting methods give symplectic integrators of arbitrary order (see \cite{McQu2002} and references therein).
\end{enumerate}
The first two techniques are based on special types of symplectic manifolds, whereas the last technique is based on a special form of the Hamiltonian function.

Recall that \emph{Poisson structures} are generalisations of symplectic structures.
The dual of a Lie algebra is called a \emph{Lie--Poisson manifold}.
It has a Poisson structure and is therefore foliated in symplectic leaves~\cite{We1983}.
Such leaves are given by \emph{coadjoint orbits}~\cite[\S\!~14]{MaRa1999}.
A \emph{Lie--Poisson integrator} on a Lie--Poisson manifold is a numerical integrator that preserves the Poisson structure and the coadjoint orbits.
When restricted to a coadjoint orbit, a Lie--Poisson integrator is thus a symplectic integrator.
The \Rattle\ method can be used to construct Lie--Poisson integrators for general Hamiltonian systems on a Lie--Poisson manifold~\cite[\S\!~VII.5]{HaLuWa2006}.

The rigid body bracket (see equation~\eqref{eq:rigid_body_bracket}) provides $\R^{3}$ with the structure of a Lie--Poisson manifold.
In this paper we develop Lie--Poisson integrators for general Hamiltonian systems on~$\R^{3}$.
The coadjoint orbits are spheres, so we also obtain symplectic integrators on~$S^{2}$.
Compared to the \Rattle\ approach, our integrators use less variables (4 instead of 10) and are less complicated since there are no constraints. 
\todo{Optimal 10 variables for Rattle: $T\SUtwo\simeq TS^3 \subset T\R^{4}$ plus two Lagrange multipliers.}

Our construction of Lie--Poisson integrators uses results on \emph{symplectic variables} (also called \emph{Clebsch variables}) and \emph{collective Hamiltonians}, found in work by \citet{MaWe1983}, \citet{GuSt1984}, and \citet{LiMa1987}.
In particular, we use \cite[Example 3.2]{MaWe1983} that the momentum map for the action of $\SUtwo$ on $\CC^{2}$ generalises the Cayley--Klein parameters for the free rigid body to arbitrary Hamiltonian problems on $\R^{3}$ (see also \cite[Exercise 15.4-1]{MaRa1999}).

The paper is organised as follows.
We first give the main result, that symplectic Runge--Kutta methods can be used to construct Lie--Poisson integrators on~$\R^{3}$.
We then describe the geometry of the problem.
Thereafter we demonstrate the main result;
we prove the core part in four different ways.
These four proofs reflect the facets of the geometry at hand.
We give a numerical example in \autoref{sec:numerical_examples} and two plausible applications in~\autoref{sec:outlook_and_applications}.
Finally, we provide two appendices detailing some computations and coordinate formulae.


\section{Main result}
\label{sec:main_results}


\noindent
The cross product provides $\R^{3}$ with the structure of a Lie algebra.
We identify the dual~$(\R^{3})^{*}$ with $\R^{3}$ by the standard inner product.
Then the Lie algebra structure on~$\R^{3}$ induces the \emph{rigid body bracket}~\cite[\S\!~10.1(c)]{MaRa1999} on~$\R^{3}$
\begin{equation}\label{eq:rigid_body_bracket}
	\poisson{F}{H}(\vect{w}) \coloneqq - \vect{w}\cdot \big(\nabla F(\vect{w})\times\nabla H(\vect{w})\big) 
.
\end{equation}
Relative to this Lie--Poisson structure, the Hamiltonian vector field corresponding to a Hamiltonian function $H\in C^{\infty}(\R^{3})$ is given by 
\begin{equation}
X_H(\vect{w}) \coloneqq \vect{w}\times \nabla H(\vect{w})
.
\end{equation}

A map $\meth[]\colon \R^{3}\to\R^{3}$ is said to \emph{preserve the Lie--Poisson structure} if
\begin{equation}\label{eq:LP_map}
	\poisson{F\circ\meth[]}{H\circ\meth[]} = \poisson{F}{H}\circ\meth[] , \qquad \forall\, F,H\in C^{\infty}(\R^{3}).
\end{equation}

Consider the map $\hopf\colon T^{*}\R^{2} \to \R^{3}$ given by
\begin{equation}\label{eq:extended_hopf}
	\hopf\colon  (q_1,q_2,p_1,p_2) \longmapsto \frac{1}{4}\begin{pmatrix}
	2q_1q_2 + 2p_1p_2 \\
	2q_1p_2 - 2q_2p_1 \\
	q_1^2 + p_1^2 - q_2^2  - p_2^2
	\end{pmatrix}
	.
\end{equation}
A Hamiltonian $H\in C^{\infty}(\R^{3})$ can be lifted to a \emph{collective Hamiltonian}
\begin{equation}\label{eq:lifted_Ham}
\pback{\hopf}{H} \coloneqq H \circ \hopf \; \in C^{\infty}(T^{*}\R^{2})
\end{equation}
on the symplectic vector space $T^{*}\R^{2}$.
By construction, the collective Hamiltonian~\eqref{eq:lifted_Ham} is constant on each \emph{fibre} $\hopf^{-1}(\set{\vect{w}})$.
The flow of the Hamiltonian vector field $X_{\pback{\hopf}{H}}$ \emph{descends} to the flow of $X_H$.
That is, this diagram commutes:
\begin{center}
\begin{tikzpicture}
	\matrix (m) [commdiag]
	{
		T^*\R^2 \& T^*\R^2  \\
		\R^{3} \& \R^{3} \\
	};
	\path[->] (m-1-1) edge node[auto]{$\exp(X_{\pback{\hopf}{H}})$} (m-1-2)
	(m-2-1) edge node[auto]{$\exp(X_H)$} (m-2-2)
	(m-1-1) edge node[auto]{$\hopf$} (m-2-1)
	(m-1-2) edge node[auto]{$\hopf$} (m-2-2)
	;
\end{tikzpicture}
\end{center}
This follows from $X_{H}(\hopf(\vect{z})) = T \hopf (X_{\pback{\hopf}{H}}(\vect{z}))$, which in turn follows from equation~\eqref{eq:hamiltonian_vf_property} and $\hopf$ being a Poisson map.

In general, a map $\Meth[]\colon T^{*}\R^{2} \to T^{*}\R^{2}$ descends if there is a map $\meth[]\colon \R^{3}\to\R^{3}$ such that this diagram commutes:
\begin{center}
\begin{tikzpicture}
	\matrix (m) [commdiag, column sep=3em]
	{
		T^*\R^2 \& T^*\R^2  \\
		\R^{3} \& \R^{3} \\
	};
	\path[->] (m-1-1) edge node[auto]{$\Meth[]$} (m-1-2)
	(m-2-1) edge node[auto]{$\meth[]$} (m-2-2)
	(m-1-1) edge node[auto]{$\hopf$} (m-2-1)
	(m-1-2) edge node[auto]{$\hopf$} (m-2-2)
	;
\end{tikzpicture}
\end{center}
Notice that $\Meth[]$ descends if and only if it maps fibres to fibres.





Our strategy for Lie--Poisson integrators on $\R^{3}$ is to integrate $X_{\pback{\hopf}{H}}$ with a symplectic Runge--Kutta method
\begin{equation}
	\Meth(X_{\pback{\hopf}{H}})\colon T^{*}\R^{2}\to T^{*}\R^{2}
\end{equation}
and then use $\pi$ to project the result back to $\R^{3}$.
To retain the structural properties of the exact flow $\exp(X_{\pback{\hopf}{H}})$, the
aim is that $\Meth(X_{\pback{\hopf}{H}})$ should descend to a Poisson integrator $\meth(H)\colon \R^{3}\to\R^{3}$.
Our main result is that any symplectic Runge--Kutta method will do.

\begin{theorem}\label{thm:main}
	Let $H\in C^{\infty}(\R^{3})$, and let $\pback{\hopf}{H}$, with~$\hopf$ given in equation~\eqref{eq:extended_hopf}, be the corresponding collective Hamiltonian on $T^{*}\R^{2}$.
	Then any symplectic Runge--Kutta method $\Meth$ applied to the Hamitonian vector field $X_{\pback{\hopf}{H}}$  descends to an integrator $\meth(H)$ on $\R^{3}$  which preserves the Lie--Poisson structure and the coadjoint orbits.
	The integrator $\meth(H)$ is consistent with $X_H$ and has the same order of convergence as $\Meth$.
	Moreover, $\meth$ is equivariant with respect to the rotation group $\SO(3)$:
	\begin{equation}
		\meth(g^{*}H) = g^{-1}\circ\meth(H)\circ g
	\end{equation}
	for all $g(\vect{w}) \coloneqq A\vect{w}$ with $A\in \SO(3)$.
\end{theorem}

So, in order to get a Lie--Poisson integrator for~$X_H$, we lift the initial data $\vect{w}_0 \in \R^{3}$ to any point $\vect{z}_0 \in \hopf^{-1}(\set{\vect{w}_{0}})$, then we integrate with a symplectic Runge--Kutta method $\Meth(X_{\pback{\hopf}{H}})$ to obtain $\vect{z}_{k+1} = \Meth(X_{\pback{\hopf}{H}})(\vect{z}_k)$, then we project to get a discrete trajectory $\vect{w}_k \coloneqq \hopf(\vect{z}_k)$ in $\R^{3}$.
\autoref{thm:main} ensures that the map $\meth(H)\colon\vect{w}_k\mapsto \vect{w}_{k+1}$ is well defined and is a Lie--Poisson integrator. 

\begin{remark}
	The Runge--Kutta method $\Meth$ in \autoref{thm:main} requires the vector field $X_{\pback{\hopf}{H}}$.
	An explicit expression for $X_{\pback{\hopf}{H}}$ is given in \autoref{sec:lifted_vf}.	
\end{remark}

\begin{remark}
	By analogy with the collective Hamiltonian $\pback{\hopf}{H}$, we call $\Meth(X_{\pback{\hopf}{H}})$ a \emph{collective integrator}.
\end{remark}

The coadjoint orbits are given by 2-spheres in $\R^{3}$ with the standard volume form as symplectic form.
Thus, the Lie--Poisson integrators in \autoref{thm:main} also yield volume preserving integrators on $S^{2}$:

\begin{corollary}\label{cor:S2integrator}
	Let $\bar H\in C^{\infty}(S^{2})$ and let $\incl\colon S^{2}\hookrightarrow \R^{3}$ be the natural inclusion.
	Let $H\in C^{\infty}(\R^{3})$ be an extension of $\bar H$, i.e., $\incl^{*}H = \bar H$.
	Then $\meth(H)$ in \autoref{thm:main} induces a volume preserving integrator $\bar\meth(\bar H)$ for $X_{\bar H}$ by $\incl\circ\bar\meth(\bar H) = \meth(H)\circ\incl$.
\end{corollary}


\section{Preliminaries}


Before the proofs we recall some facts about Poisson manifolds, group actions, and group invariant vector fields.
This will make the geometry of the problem more transparent.

\subsection{Poisson manifolds} 
\label{par:poisson_manifolds}

If $(P_1,\{\cdot,\cdot\}_1)$ and $(P_2,\{\cdot,\cdot\}_2)$ are two Poisson manifolds~\cite[\S\!~10]{MaRa1999}, then a map $R\colon P_1\to P_2$ is a \emph{Poisson map} if
\begin{equation}\label{eq:poisson_map}
	\{ F\circ R,G\circ R \}_{1} = \{ F,G \}_{2}\circ R, \qquad\forall F,G\in C^{\infty}(P_2).
\end{equation}
Furthermore, if $H\in C^{\infty}(P_2)$ then
\begin{equation}\label{eq:hamiltonian_vf_property}
	T R\circ X_{H\circ R} = X_H\circ R .
\end{equation}
See, e.g.,\ \citet[Prop.~10.3.2]{MaRa1999}.


\subsection{Group actions on $\CC^{2}$} 
\label{par:group_action_on_}

Since $\CC^{2}$ is a complex Hilbert space it is equipped with a symplectic structure (see, e.g.,~\citet[\S\!~5.2]{MaRa1999}).
Under the vector space isomorphism
\begin{equation}\label{eq:identification_C2TstarR2}
	T^{*}\R^{2} \ni (q_0,q_1,p_0,p_1) \longmapsto  (q_0 + \ii p_0,q_1 + \ii p_1) \in \CC^{2}.
\end{equation}
this symplectic structure is exactly the canonical symplectic structure of $T^{*}\R^{2}$.

The group $\SUtwo$ of complex special unitary $2\times 2$~matrices acts on $\CC^2$ by matrix multiplication.
Since this action preserves the Hermitian inner product on $\CC^{2}$ it is automatically symplectic.
Recall that the Lie algebra $\sutwo$ of $\SUtwo$ consists of skew Hermitian and trace free complex matrices.
The cross product algebra $(\R^{3},\times)$ can be identified with $\su(2)$ by the Lie algebra isomorphism~\cite[\S\!~9, p.~302]{MaRa1999}
\begin{equation}\label{eq:su2R3isomorphism}
	\R^{3}\ni (x_1,x_2,x_3) \longmapsto \frac{1}{2}\begin{pmatrix}
		-\ii x_3 & -\ii x_1 -x_2 \\
		-\ii x_1 + x_2 & \ii x_3
	\end{pmatrix} \in \su(2).
\end{equation}
By the isomorphisms~\eqref{eq:identification_C2TstarR2}, \eqref{eq:su2R3isomorphism}, and $(\R^{3})^{*}\simeq\R^{3}$, the projection map~\eqref{eq:extended_hopf} can be interpreted as a map $\hopf\colon\CC^{2}\to\sutwo^{*}$.
If we take this point of view, then $\hopf$ is the \emph{momentum map} associated with the action of $\SUtwo$ on $\CC^{2}$ (see \autoref{sec:momentum_map_computed} for details).
This momentum map is $\SUtwo$-equivariant, i.e., $\hopf(A\cdot \vect{z}) = A\cdot\hopf(\vect{z})$, where $A \in \SUtwo$ acts on $\sutwo^{*}$ by the coadjoint action $\Ad_{A^{-1}}^{*}\colon\sutwo^{*}\to\sutwo^{*}$ (this action corresponds to the action of $\SO(3)$ on $\R^{3}$~\cite[Prop.~9.2.19]{MaRa1999}).
A consequence is that $\hopf$ is a Poisson map~\cite[Th.~12.4.1]{MaRa1999}.

\begin{remark}
	The map $\hopf$ takes spheres in $\CC^{2}$ of radius $r$ into spheres in $\R^{3}\simeq\su(2)^{*}$ of radius $r^{2}/4$.
	That is:
	\begin{equation}\label{eq:hopf_map_sphere_radius}
		\norm{\hopf(\vect{z})}_{2} = \frac{1}{4}\norm{\vect{z}}_{2}^{2}.
	\end{equation}
	Much related to $\pi$ is the map
	\begin{equation}\label{eq:original_hopf_map}
		\CC^{2} \supset S^{3} \ni \vect{z} \longmapsto 4\hopf(\vect{z})\in S^{2}\subset \R^{3},
	\end{equation}
	which a Riemannian submersion.
	This map yields the classical \emph{Hopf fibration}, i.e., the fibration of $S^{3}$ in circles over $S^2$.
\end{remark}

The Lie group of unitary complex numbers $\Uone$ acts on $\CC^{2}$ by 
\begin{equation}\label{eq:U1group_action}
	(z_1, z_2)\cdot\ee^{\ii \theta}  = (z_1 \ee^{\ii \theta}, z_2 \ee^{\ii \theta})
	.
\end{equation}
We denote this action by $\Uoneact \colon \CC^2 \to \CC^2$.
It preserves each fibre $\hopf^{-1}(\set{\vect{w}})$.
That is, 
\begin{equation}
\label{eq:fibre_action}
\hopf\circ\Uoneact=\hopf
.
\end{equation}
In fact, the action is transitive in the fibres, so $\Uone$ gives a parametrisation of each fibre.
In particular, if $\Meth[]$ is $\Uone$-equivariant, i.e., $\Meth[](\vect z\cdot\ee^{\ii\theta}) = \Meth[](\vect{z})\cdot\ee^{\ii\theta}$, then $\Meth$ is descending.



The $\Uone$~action preserves the Hermitian inner product, so it is symplectic.
The orbits are generated by the flow of a Hamiltonian vector field $X_M$, with Hamiltonian function 
\begin{equation}
\label{eq:momentum}
M(z_1,z_2) = \abs{z_1}^2 + \abs{z_2}^2
.
\end{equation}
Thus, $M$ is the momentum map for the $\Uone$~action and may be regarded as a map $M\colon \CC^{2}\to \uone^{*}\simeq \R$.

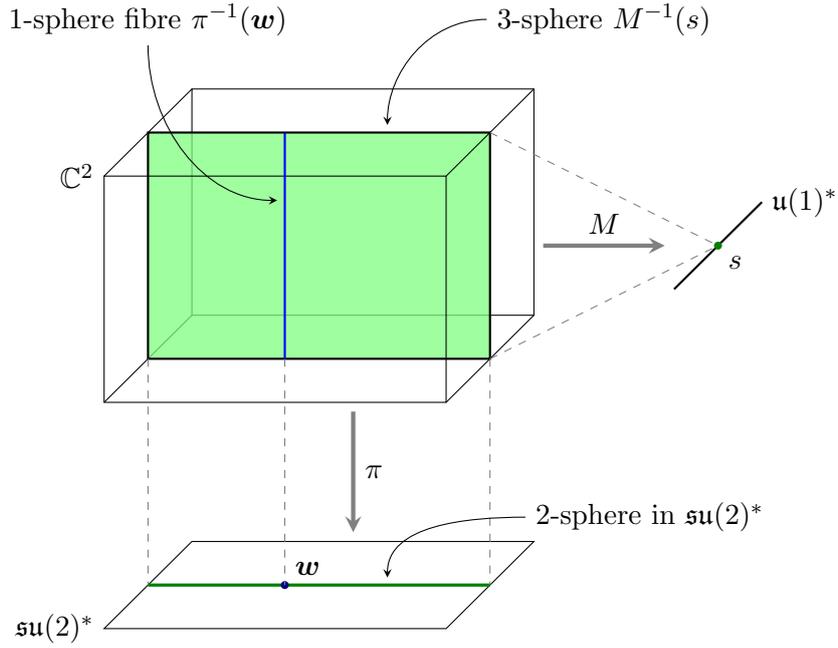
\begin{figure}
	\centering
		
\newcommand*\dx{1.5}
\newcommand*\dy{-1}
\newcommand*\upos{2.5,.5}
\newcommand*\drawdepthline[2][]{
	\draw[#1] (A#20) -- (A#21);
}
\newcommand*\drawdepth[2]{
{
\path[draw,#2] (A00#1) -- (A01#1) -- (A11#1) -- (A10#1) -- cycle;
}
}
\begin{tikzpicture}[scale=3,
>=stealth,
projline/.style={gray, dashed},
projarrow/.style={->, ultra thick, gray, shorten <=.7cm, shorten >=.7cm,},
labelarrow/.style={->, shorten >=.1cm},
labeltext/.style={shape=rectangle,},
]
\coordinate (A000) at (0,0,0);
\coordinate (A001) at (0,0,1);
\coordinate (A010) at (0,1,0);
\coordinate (A011) at (0,1,1);
\coordinate (A100) at (\dx,0,0);
\coordinate (A101) at (\dx,0,1);
\coordinate (A110) at (\dx,1,0);
\coordinate (A111) at (\dx,1,1);
\coordinate (A00m) at (0,0,.5);
\coordinate (A01m) at (0,1,.5);
\coordinate (A10m) at (\dx,0,.5);
\coordinate (A11m) at (\dx,1,.5);
\coordinate (Am0) at (\upos,0);
\coordinate (Am1) at (\upos,1);
\coordinate (Am) at (\upos,.5);
\coordinate (A0d0) at (0,\dy,0);
\coordinate (A0d1) at (0,\dy,1);
\coordinate (A1d0) at (\dx,\dy,0);
\coordinate (A1d1) at (\dx,\dy,1);
\coordinate (A0dm) at (0,\dy,.5);
\coordinate (A1dm) at (\dx,\dy,.5);
\coordinate (Ap0m) at ($(A00m)!.4!(A10m)$);
\coordinate (Ap1m) at ($(A01m)!.4!(A11m)$);
\coordinate (Apdm) at ($(A0dm)!.4!(A1dm)$);
	\drawdepth{0}{}
	\drawdepthline{10}
	\drawdepthline{00}
	\drawdepthline{01}
	\drawdepth{m}{thick,fill=green!50!white,fill opacity=.75}
	\drawdepthline{11}
	\drawdepth{1}{}
	\draw (A0d0) -- (A0d1) -- (A1d1) -- (A1d0) -- cycle;
	\draw[very thick,green!50!black] (A0dm) -- (A1dm);
	\drawdepthline[thick]{m}
	\fill[green!50!black] (Am) circle[radius=.5pt] ;
	\fill[blue!50!black] (Apdm) circle[radius=.5pt] ;
	\draw[projline] (Ap0m) -- (Apdm);
	\draw[projline] (A00m) -- (A0dm);
	\draw[projline] (A10m) -- (A1dm);
	\draw[projline] (A10m) -- (Am);
	\draw[projline] (A11m) -- (Am);

	\draw[thick,blue] (Ap0m) -- (Ap1m);

	\draw[projarrow] ($(A10m)!.5!(A11m)$) -- node[auto,black]{$M$} (Am);
	\draw[projarrow] ($(A00m)!.6!(A10m)$) -- node[auto,black]{$\hopf$} ($(A0dm)!.6!(A1dm)$);

	\node[left] at (A011) {$\CC^2$};
	\node[left] at (A0d1) {$\su(2)^*$};
	\node[right] at (Am0) {$\uone^*$};
	\node[below right] at (Am) {$s$};

	\node[labeltext] (twospherelabel) at ($(A1dm)+(.7,.3,0)$) {2-sphere in $\su(2)^*$};
	\draw[labelarrow]   (twospherelabel)  to[in=90,out=180] ($(A0dm)!.7!(A1dm)$);

	\node[labeltext] (threespherelabel) at ($(A11m)+(.5,.5,0)$) {3-sphere $M^{-1}(s)$};
	\draw[labelarrow] (threespherelabel) to[out=-180,in=90] ($(Ap1m)!.5!(A11m)$);

	\node[labeltext] (fibrelabel) at ($(A11m)+(-1.5,.5,0)$) {1-sphere fibre $\hopf^{-1}(\vect{w})$};
	\draw[labelarrow] (fibrelabel) to[out=-90, in=180] ($(Ap1m)!.3!(Ap0m)$);

	\node[above right] at (Apdm) {$\vect{w}$};


\end{tikzpicture}
	\caption[Hopf]{
		Visualisation of the relation between the spaces involved in the proof.
		Recall that $\hopf$ is defined by \eqref{eq:extended_hopf} and the identification \eqref{eq:su2R3isomorphism}, and that the momentum map $M$ is defined by \eqref{eq:momentum}.
	}
	\label{fig_hopf}
\end{figure}

An illustrative summary of the spaces involved is provided in \autoref{fig_hopf}.

\subsection{Invariant vector fields}

The pullback by a diffeomorphism $\Psi$ of a vector field $X$ is defined as 
\begin{equation}\label{eq:pullback_vf}
	\Psi^{*}X := T\Psi^{-1}\circ X\circ \Psi.	
\end{equation}
The following result on arbitrary Runge--Kutta methods is useful:
\begin{lemma}
\label{lem:descendingRK}
	Let $G$ be a Lie group which acts linearly on $\CC^{n}$ by an action map $\Psi_g$.
	Let $X$ be a $G$-invariant vector field on $\CC^{n}$, that is:
	\begin{equation}
		\Psi_g^{*}X = X, \quad \forall\, g\in G.
	\end{equation}
	Further, let $\Meth(X)\colon\CC^{n}\to\CC^{n}$ be a Runge--Kutta method applied to~$X$.
	Then $\Meth(X)$ is $G$-invariant, that is:
	\begin{equation}
		\Psi_g^{-1}\circ\Meth(X)\circ\Psi_g = \Meth(X), \quad \forall\, g\in G.
	\end{equation}
\end{lemma}
\begin{proof}
	Runge--Kutta methods are equivariant with respect to affine transformations.
	Since $\Psi_g$ is linear this implies
	\begin{equation}
		\Psi_g^{-1}\circ\Meth(X)\circ\Psi_g = \Meth(\Psi_g^{*}X).
	\end{equation}
	The vector field $X$ being $G$-invariant thus implies
	\begin{equation}
		\Psi_g^{-1}\circ\Meth(X)\circ\Psi_g = \Meth(X),		
	\end{equation}
	i.e., $\Meth(X)$ is $G$-invariant.
\end{proof}


Notice that that if $H \in C^{\infty}(\su(2)^*)$ then the Hamiltonian vector field $X_{\pback{\hopf}{H}}$ on $\CC^{2}$ is $\Uone$-invariant.
Indeed, the action $\Uoneact$ is symplectic, so equation~\eqref{eq:hamiltonian_vf_property} and equation~\eqref{eq:fibre_action} yield
\begin{equation}
\label{eq:lifted_vf_invariant}
	\Uoneact^{*}X_{\pback{\hopf}{H}} = X_{\pback{\Uoneact}{\pback{\hopf}{H}}} = X_{\pback{\hopf}{H}}.
\end{equation}

\section{Proofs}

\autoref{thm:main} is proved in three steps:
\begin{enumerate}
	\item In \autoref{sub:proof_descend} we prove, in four different ways, that $\Meth(X_{\pback{\hopf}{H}})$ descends to a Poisson map $\meth(H)$ which preserves coadjoint orbits.
	
	\item In \autoref{sub:convergence_order} we prove that $\meth(H)$ has the same convergence order as $\Meth(X_{\pback{\hopf}{H}})$.
	\item In \autoref{sub:equivariance} we prove that $\meth(H)$ is $\SUtwo$-equivariant.
\end{enumerate}


\subsection{Descend and coadjoint orbits} 
\label{sub:proof_descend}

\begin{proposition}\label{prop:descend}
	The map $\Meth(X_{\pback{\hopf}{H}})\colon\CC^{2}\to\CC^{2}$ descends to a map $\meth(H)\colon\su(2)^{*}\to\su(2)^{*}$ which preserves the Lie--Poisson structure and the coadjoint orbits.
\end{proposition}
We give four different proofs of this result.
Each proof suggests its own generalisation to other Lie--Poisson manifolds.
Such generalisations are explored in a forthcoming paper~\cite{McMoVe2013b}.

To start with, notice the special coadjoint orbit at the origin;
it has dimension zero, whereas all the others have dimension two.
Furthermore, the fibre above the origin in $\su(2)^{*}$ is the origin in $\CC^{2}$, so this \emph{trivial fibre} has dimension zero, whereas all the other fibres have dimension one.
The origin in $\su(2)^{*}$ is an equilibrium point for every Hamiltonian vector field $X_H$ on $\su(2)^{*}$, so the origin in $\CC^{2}$ is an equilibria point for every lifted Hamiltonian vector field $X_{\pback{\hopf}{H}}$ on $\CC^{2}$.
Since every Runge--Kutta method exactly preserves equilibria points, the result in \autoref{prop:descend} is true at the origin.
Thus, it remains to prove the result on $\CC^{2}\setminus \set{0}$.
Removing the origin simplifies the proofs, since $\pi\colon\CC^{2}\setminus \set{0} \to \su(2)^{*}\setminus \set{0}$ is a submersion, so if $\Meth[]\colon\CC^{2}\setminus\set{0}\to\CC^{2}\setminus\set{0}$ descends to $\meth[]\colon\su(2)^{*}\setminus\set{0}\to\su(2)^{*}\setminus\set{0}$ and $\Meth[]$ is smooth, then $\meth[]$ is also smooth (and correspondingly for smooth functions).

The following result is used in the first, second, and fourth proof.

\begin{lemma}
\label{lma:momentum}
Let $\Meth$ be a symplectic Runge--Kutta method on $\CC^{2}$.
Then the momentum map $M$ in equation~\eqref{eq:momentum} is conserved by $\Meth(X_{\pback{\hopf}{H}})$.	
\end{lemma}
\begin{proof}
	The momentum $M$ associated with the $\Uone$~action is quadratic.
	Since $M$ is an invariant of $X_{\pback{\hopf}{H}}$, and since symplectic Runge--Kutta methods conserve quadratic invariants~\cite[Th.~VI.7.6]{HaLuWa2006}, we get that $\Meth(X_{\pback{\hopf}{H}})$ conserves~$M$.
\end{proof}

The following result is used in the first and fourth proof.

\begin{lemma}
\label{lma:sutwo_equivariance}
Let $\Meth[]\colon \CC^{2}\to\CC^{2}$ be a map that descends and conserves the momentum map~$M$ in equation~\eqref{eq:momentum}.
Then the descending map $\meth[]\colon \su(2)^{*}\to\su(2)^{*}$ preserves the coadjoint orbits.
\end{lemma}

\begin{proof}
The action of $\SUtwo$ on $\CC^{2}$ is norm preserving.
For every $A\in\SUtwo$, $\vect{z}$ and $A\cdot\vect{z}$ belong to the same level set of $M$ because $M(\vect z) = \norm{\vect z}^{2}$.


Let $\vect{z}' \coloneqq \meth[](\vect{z})$.
Then $\vect{z}$ and $\vect{z}'$ belong to the same level set of $M$.
Since $\SUtwo$ acts transitively on that level set~\cite[Eq.~9.2.15]{MaRa1999}, there is an $A \in \SUtwo$ such that $\vect{z}' = A\cdot \vect{z}$.
By equivariance of the momentum map $\hopf$, we conclude that $\hopf(\vect{z}') = A\cdot \hopf(\vect{z})$, so the descending map $\meth[]$ preserves the coadjoint orbit.
\end{proof}

\paragraph{Symmetry + conservation of momentum} 
\label{sub:proof_symmetry_momentum}

The first proof we consider is based on two key properties: (i) the Hamiltonian $\pback{\hopf}{H}$ is invariant with respect to the $\Uone$~action, and (ii) symplectic Runge--Kutta methods conserve the momentum~$M$ associated with the $\Uone$~action.

\begin{proof}[First proof of \autoref{prop:descend}]

	We observed in equation~\eqref{eq:lifted_vf_invariant} that $X_{\pback{\hopf}{H}}$ is $\Uone$-equivariant.
	Thus, it follows from \autoref{lem:descendingRK} that $\Meth(X_{\pback{\hopf}{H}})$ is $\Uone$-equivariant.
	The $\Uone$~action parameterises the fibres, which implies that $\Meth(X_{\pback{\hopf}{H}})$ descends to a map $\meth( H)$.
	This map preserves the Lie--Poisson structure, since $\Meth(X_{\pback{\hopf}{H}})$ is symplectic and $\hopf$ is a Poisson map.

	Finally, we conclude that the descending method $\meth(H)$ preserves the coadjoint orbits by appealing to \autoref{lma:momentum} and \autoref{lma:sutwo_equivariance}.
\end{proof}


\paragraph{Backward error analysis + conservation of momentum} 
\label{sub:proof_BEA_momentum}

The idea in this proof is to use \emph{backward error analysis} of numerical integrators.
More precisely, we follow the framework developed in~\cite[\S\!~IX]{HaLuWa2006}.
Throughout this section and the next, a truncated modified Hamiltonian for the symplectic Runge--Kutta method $\Meth(X_{\pback{\hopf}{H}})$ is denoted $\modifham$ (i.e., an arbitrary truncation of the series in~\cite[Eq.~IX.3.4]{HaLuWa2006}).

\begin{proof}[Second proof of \autoref{prop:descend}]
	In terms of backward error analysis, the result of \autoref{lma:momentum} is that
\begin{equation}
	\pair{\ud M}{X_{\modifham}} = 0
	.
\end{equation}
We now have
\begin{equation}
	0 = \pair{\ud M}{X_{\modifham}} = \{M,\modifham \} = -\{\modifham,M \} = -\pair{\ud\modifham}{X_M}.
\end{equation}
Since $X_M$ generates the $\Uone$~orbits, and since the $\Uone$~orbits parameterises the fibres, the equality $\pair{\ud\modifham}{X_M} = 0$ implies that $\modifham$ is constant on the fibres.
Thus, $\modifham =  \pback{\hopf}{\tilde{H}_h}$ for some function $\tilde{H}_h$ on $\su(2)^{*}$.
When restricted to $\su(2)^{*}\setminus\set{0}$ this function is smooth.
Therefore, by equation~\eqref{eq:hamiltonian_vf_property}, the Hamiltonian vector field $X_{\modifham}$ restricted to $\CC\setminus\set{0}$ descends to the Hamiltonian vector field $X_{\tilde{H}_h}$ on $\su(2)^{*}\setminus\set{0}$, which implies that the flow $\exp(h X_{\modifham})$ descends to $\exp(h X_{\tilde{H}_h})$.
Since this is true for all truncations of modified Hamiltonians $\modifham$, it follows from \cite[Th.~IX.5.2]{HaLuWa2006} that $\Meth(X_{\pback{\hopf}{H}})$ descends to a Poisson map which preserves coadjoint orbits.
\end{proof}


\paragraph{Backward error analysis + trivial cohomology} 
\label{sub:proof_BEA_cohomology}

In this proof we again use backward error analysis, but instead of conservation of momentum we use equivariance of Runge--Kutta methods and that a certain cohomology class is trivial.
The notation is the same as in the previous proof.

\begin{proof}[Third proof of \autoref{prop:descend}]

	To conduct backward error analysis we examine the group of $\Uone$-equivariant diffeomorphism and its formal Lie algebra.
	The set of all $\Uone$-equivariant diffeomorphisms on $\CC^{2}$ forms a group denoted~$\Diff_{\Uone}(\CC^{2})$.
	The corresponding algebra is given by the space of $\Uone$-invariant vector fields and is denoted $\Xcal_{\Uone}(\CC^{2})$.
	Since $X_M$ generates the $\Uone$~orbits, it follows that 
	\begin{equation}
		X\in \Xcal_{\Uone}(\CC^{2}) \iff [X,X_M] = 0.
	\end{equation}

	 We have already seen in equation~\eqref{eq:lifted_vf_invariant} that $X_{\pback{\hopf}{H}} \in \Xcal_{\Uone}(\CC^{2})$.
	Since Runge--Kutta methods are equivariant with respect to affine transformations, and since the $\Uone$~action is linear, \autoref{lem:descendingRK} gives that $\Meth(X_{\pback{\hopf}{H}}) \in \Diff_{\Uone}(\CC^{2})$.
	In terms of backward error analysis, this means that $X_{\modifham}\in\Xcal_{\Uone}(\CC^{2})$, i.e.,
	\begin{equation}
		[X_{\modifham},X_{M}] = 0.
	\end{equation}
	We now get
	\begin{equation}
		 X_{\{M, \modifham  \}} = [X_{\modifham},X_{M}] = 0
		,
	\end{equation}
	which implies that $\{M, \modifham\}$ is constant, so $\pair{\dd \modifham}{X_{M}} = c$ for some constant $c\in\R$.

	In order to show that~$c=0$, we consider one fixed $\Uone$-orbit $O$, for instance the one passing through $\vect{z}_0 = (1,0)$.
	Let $\incl\colon O\hookrightarrow \CC^{2}$ be the inclusion of $O$ into $\CC^2$.
	Since $O$ is a closed manifold we get
	\begin{equation}
	\label{eq:cohomology}
		\begin{split}
			0 &= \int_{O} \ud \incl^{*}\modifham \\
			& = \int_{O} \ud \Big( \modifham\big(\Uoneact(\vect{z}_0)\big) \Big)\\
			& = \int_{0}^{2\pi} \pair{\ud \modifham\big(\Uoneact(\vect{z}_0) \big)}{\frac{\ud \Uoneact(\vect{z}_0)}{\ud\theta}}\udd\theta \\
			& = \int_{0}^{2\pi} \underbrace{\pair{\ud \modifham\big(\Uoneact(\vect{z}_0) \big)}{X_{M}\big(\Uoneact(\vect{z}_0)\big)}}_{c}\udd\theta \\
			&= c \int_{0}^{2\pi}\ud\theta = 2 \pi c , 
		\end{split}
	\end{equation}
	which implies $c=0$.
	Thus, $\pair{\dd \modifham}{X_{M}} = 0$, so $\modifham$ is constant on the fibres.
	Repeating the last part of the previous proof, $\Meth(X_{\pback{\hopf}{H}})$ descends to a Poisson map which preserves coadjoint orbits.
\end{proof}

The last part of the proof hinges on a property of the fibration related to the cohomology of the corresponding quotient space.
That property is that if a one-form is a pull-back by $\hopf$ and is closed, then it is the exterior differential of a function which is also a pullback by $\hopf$.
The meaning of equation~\eqref{eq:cohomology} is to show that this holds, in other words, that the cohomology associated to the fibration is trivial.

\begin{remark}
	Concerning the cohomology above, we believe that there is a minor mistake in \cite[Th.~IX.5.7]{HaLuWa2006}. 
	That theorem states that any left invariant symplectic integrator on the cotangent bundle $T^*G$ of a Lie group $G$ descends to a Poisson integrator on $\mathfrak{g}^{*}$ which preserves coadjoint orbits.
	The statement that coadjoint orbits are preserved is true only if the cohomology class of $G$-invariant one-forms on $G$ is trivial.
	This is not always the case, not even locally.
\end{remark}


\paragraph{Collective approach} 
\label{sub:proof_collective}

We use that $\SUtwo$ acts symplectically on $\CC^2$ and that the Hopf map is the equivariant momentum map for the action.
This proof does not use backward error analysis.

\begin{proof}[Fourth proof of \autoref{prop:descend}]

We know from \autoref{lma:momentum} that $M$ is preserved.
Since $\Meth$ is symplectic, it preserves the symplectic orthogonal directions.
Indeed, suppose that $X$ and $Y$ are two symplectically orthogonal vectors, i.e., $\omega(X,Y) = 0$.
By symplecticity of $\Meth$ we have 
\begin{equation}
	\omega(T\Meth X, T\Meth Y) = \omega(X,Y) = 0,
\end{equation}
so $T\Meth(X)$ and $T\Meth(Y)$ are symplectic orthogonal.

The symplectic orthogonal direction to a tangent space of a level set of $M$ is the fibre direction at that point, since any tangent vector $X$ of a level set of $M$ fulfils $\omega(X_M, X) = \pair{\ud M}{X} = 0$.
If $\vect{z}\in\CC\setminus\set{0}$ then, by a dimension argument, $X_M(\vect{z})\neq 0$ is the only direction which is symplectic orthogonal to the tangent space at $\vect{z}$ of the level sets of $M$.

Take two points $\vect{z}$, $\vect{z}'$ in a non-trivial fibre.
There is a curve in that fibre which joins these two points.
The image of that curve by $\Meth$ is another curve, which is tangent at all its points to the fibre direction.
As a result, that image curve also lies in a fibre.
Its end points $\Meth(\vect{z})$ and $\Meth(\vect{z}')$ are thus in the same fibre, from which we conclude that $\Meth$ is fibre preserving.

Finally, we conclude that the descending method $\meth$ preserves the coadjoint orbits by appealing to \autoref{lma:momentum} and \autoref{lma:sutwo_equivariance}.
\end{proof}



\subsection{Convergence order} 
\label{sub:convergence_order}

\begin{proposition}\label{prop:order}
	Let $\Meth(X_{\pback{\hopf}{H}})$ be an integrator of order~$k$ for $X_{\pback{\hopf}{H}}$ that descends to~$\meth(H)$.
	Then $\meth(H)$ is an integrator of order~$k$ for $X_H$.
\end{proposition}

\begin{proof}
	We have
	\begin{equation}
		\Meth(X_{\pback{\hopf}{H}})(\vect z) = \exp(h X_{\pback{\hopf}{H}})(\vect{z}) + \mathcal{O}(h^{k+1}).
	\end{equation}
	Since $\hopf$ is smooth this yields
	\begin{equation}
		\big(\hopf\circ\Meth(X_{\pback{\hopf}{H}})\big)(\vect{z}) = \big(\pi\circ\exp(h X_{\pback{\hopf}{H}})\big)(\vect{z})  + \mathcal{O}(h^{k+1}).
	\end{equation}
	Next, since $\Meth(X_{\pback{\hopf}{H}})$ and $\exp(h X_{\pback{\hopf}{H}})$ descend we get
	\begin{equation}
		\meth(H)(\hopf(\vect{z})) = \exp(X_{H})(\hopf(\vect{z})) + \mathcal{O}(h^{k+1}).
	\end{equation}
	The result now follows since $\pi$ is surjective.
\end{proof}


\subsection{Equivariance} 
\label{sub:equivariance}

Let $\Psi_A$, $\psi_A$ denote the action of $A\in\SUtwo$ on $\CC^{2}$, $\su(2)^{*}$ respectively (see \autoref{par:group_action_on_}).


\begin{proposition}

Assume that $\Meth(X_{\pback{\hopf}{H}})$ descends to $\meth(H)$.
Further, assume that $\Meth$ is $\SUtwo$-equivariant, i.e.,
\begin{equation}
	\Meth(\Psi_A^{*}X) = \Psi_A^{-1}\circ\Meth(X)\circ \Psi_A,
\end{equation}
for all vector fields $X$ on $\CC^{2}$.
Then $\meth(H)$ is also $\SUtwo$-equivariant, i.e.,
\begin{equation}
	\meth(\pback{\psi_A}{H}) = \psi_A^{-1}\circ\meth(H)\circ \psi_A .
\end{equation}

\end{proposition}

\begin{proof}
From equation~\eqref{eq:hamiltonian_vf_property} we get $X_{\pback{\Psi_A}{H}} = \Psi_A^{*} X_H$, because $\Psi_A$ is a Poisson map.
By $\SUtwo$-equivariance of $\hopf$ one has $\hopf \circ \Psi_A = \psi_A \circ \hopf$.
Therefore
\begin{align}
\Meth(X_{\pback{\hopf}{\pback{\psi_A}{H}}}) &= \Meth(X_{\pback{\Psi_A}{\pback{\hopf}{H}}}) &\text{by equivariance of $\hopf$}\\
&= \Meth(\Psi_A^{*}  X_{\pback{\hopf}{H}}) &\text{because $\Psi_A$ is Poisson}\\
&= \Psi_A^{-1} \circ \Meth(X_{\pback{\hopf}{H}}) \circ \Psi_A &\text{by equivariance of $\Meth$}
\end{align}
Now $\meth(H)$ is defined by $\meth(H) \circ \hopf = \hopf \circ \Meth(X_H)$ so we have
\begin{align}
\meth(\pback{\psi_A}{H}) \circ \hopf &= \hopf \circ \Meth(X_{\pback{\hopf}{\pback{\psi_A}{H}}}) \\
&= \hopf \circ \Psi_A^{-1} \circ \Meth(X_{H}) \circ \Psi_A \\
&= \psi_A^{-1}\circ (\hopf \circ \Meth(X_{H}))\circ \Psi_A \\
&= \psi_A^{-1} \circ \meth(H) \circ \hopf \circ \Psi_A \\
&= \psi_A^{-1} \circ \meth(H)\circ \psi_A \circ \hopf
\end{align}
We conclude that $\meth(\pback{\psi_A}{H}) = \psi_A^{-1}\circ \meth(H) \circ \psi_A$ because $\hopf$ is surjective.
\end{proof}


\section{Numerical example} 
\label{sec:numerical_examples}

The free rigid body is a standard test problem for Lie--Poisson integrators on~$\R^{3}$~\cite[\S\!~VII.5]{HaLuWa2006}.
The Hamiltonian is
\begin{equation}\label{eq:Ham_rigid_body}
	H(w_1,w_2,w_3) = \frac{1}{2}\sum_{i=1}^{3} \frac{w_i^2}{I_i} 
,
\end{equation}
for given moments of inertia $I_1$, $I_2$, and $I_3$.

\begin{figure}
	\centering
\subfloat[Trajectories]{\label{fig_rigid_traj}\includegraphics[width=.39\textwidth]{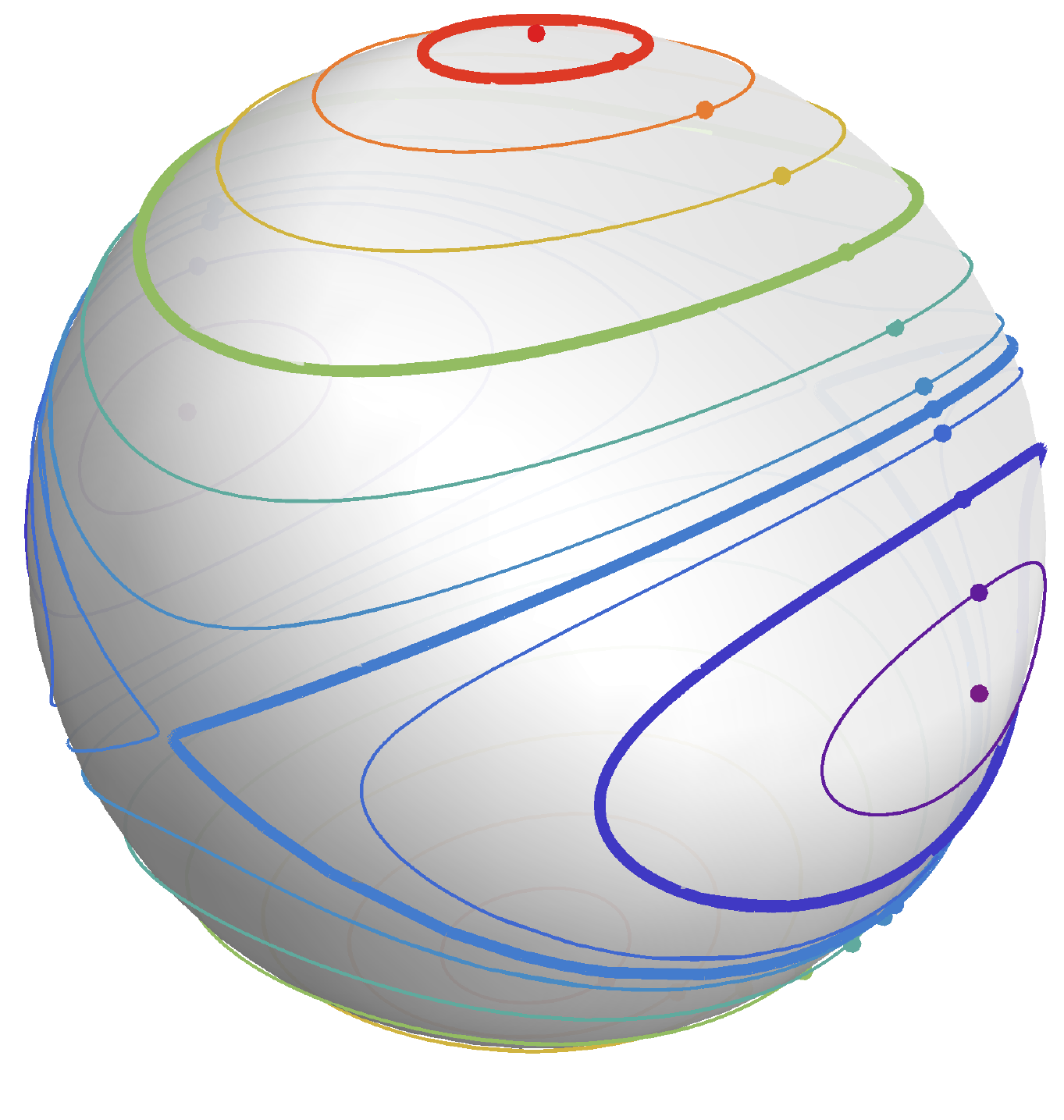}}
\subfloat[Energy errors]{\label{fig_rigid_energy}\input{fig_energy.tex}}
\caption{
Simulation of the free rigid body problem with Hamiltonian~\eqref{eq:Ham_rigid_body}.
The moments of inertia are $I_1 = 3/2$, $I_2 = 1$, and $I_3 = 1/2$.
We use the midpoint rule as collective integrator.  
The step size is $h = 5 \times 10^{-2}$.
\protect\subref{fig_rigid_traj} Some computed trajectories projected back to $\R^{3}$.
Up to machine precision the computed trajectories stay on the sphere, i.e., on a coadjoint orbit.
The computed trajectories for initial data lying on heteroclinic orbits are almost correct.
\protect\subref{fig_rigid_energy} Energy error $H(\vect{w}_k)-H(\vect{w}_0)$ for some computed solutions (corresponding to the thicker trajectories in \autoref{fig_rigid_traj}).
The error is larger for trajectories close to the heteroclinic orbits.
We observe that the energy is maximal at the initial points.
}
\label{fig_rigidbody}
\end{figure}

This problem is integrated with the collective integrator $\Meth(X_{\pback{\hopf}{H}})$, using the implicit midpoint rule as symplectic Runge--Kutta method.
We chose various initial data lifted to fibres above the sphere of radius one.
The resulting trajectories, projected back to $\R^{3}$, are plotted in \autoref{fig_rigid_traj}.
The energy evolution is plotted in \autoref{fig_rigid_energy}.

\section{Application outlook} 
\label{sec:outlook_and_applications}

As pointed out in \autoref{cor:S2integrator}, a significant feature of the collective integrators described in this paper is that they can be used for \emph{general} Hamiltonian systems on $S^2$.
In this section we propose two potential applications, both concerned with global atmospheric dynamics.
These applications will be developed in detail in future work.

\paragraph{Particle tracking for incompressible flows on $S^{2}$} 
\label{par:particle_tracking_for_incompressible_flows}

The outcome of a global weather simulation is often a time dependent divergence free vector field on $S^{2}$ describing the infinitesimal evolution of fronts.
If a smooth enough, e.g., $C^1$, interpolation of the computed vector field is known, then the collective integrator in \autoref{thm:main} can be used to accurately track the flow of individual particles, while exactly conserving mass (since the computed flow is volume preserving).


\paragraph{Point vortex dynamics} 
\label{par:point_vortex_dynamics}

The study of point vortices on two-dimensional manifolds is a ``classical mathematics playground'' according to \citet{Ha2007}.
More than a playground, point vortices on the sphere provide models for various phenomena in global geophysical flows, such as Jupiter's great red spot, or the Gulf Stream rings in the ocean.
Mathematically, the dynamics of $N$ point vortices on the sphere is described by a Hamiltonian system on $(S^{2})^{N}$ equipped with the product symplectic structure inherited from~$S^{2}$.
The Hamiltonian is non-separable so position/momentum based splitting methods can not be used.
However, by $N$ copies of the momentum map~\eqref{eq:extended_hopf}, the Hamiltonian on $(S^{2})^{N}$ can be lifted to a collective Hamiltonian on $(T^*\R^{2})^{N}\simeq T^*\R^{2N}$, which can be integrated using a symplectic Runge--Kutta method.
From \autoref{thm:main} and some book-keeping related to the product symplectic structure, it follows that the method will descend to a symplectic integrator on $(S^{2})^{N}$.




\appendix

\section{Computation of the momentum map $\pi$}\label{sec:momentum_map_computed}

In this appendix we show, by explicit computations, that the map $\pi$ in equation~\eqref{eq:extended_hopf} is the momentum map for the action of $\SUtwo$ on $\CC^{2}$.

First, the infinitesimal generator corresponding to an element $\xi\in\su(2)$ is a vector field $\xi_{\CC^{2}}$ on $\CC^{2}$ given by
\begin{equation}
	\xi_{\CC^{2}}(\vect{z}) = \xi \vect{z}
\end{equation}
By the algebra isomorphism~\eqref{eq:su2R3isomorphism}, $\xi\in\su(2)$ corresponding to $(x_1,x_2,x_3)\in\R^{3}$ is given by
\begin{equation*}
	\xi = \frac{1}{2}\begin{pmatrix}
		-\ii x_3 & -\ii x_1 -x_2 \\
		-\ii x_1 + x_2 & \ii x_3
	\end{pmatrix}
	.
\end{equation*}
The pullback of the vector field $\xi_{\CC^{2}}$ by the symplectic isomorphism~\eqref{eq:identification_C2TstarR2} between $T^*\R^{2}$ and $\CC^{2}$ yields the vector field on $T^*\R^{2}$ given by
\begin{equation}
	\xi_{T^*\R^{2}}(q_1,q_2,p_1,p_2)	= \frac{1}{2}\begin{pmatrix}
		\phantom{+}p_2 x_1 - q_2 x_2 + p_1 x_3 \\
		\phantom{+}p_1 x_1 + q_1 x_2 - p_2 x_3 \\
		-q_2 x_1 - p_2 x_2 - q_1 x_3 \\
		-q_1 x_1 + p_1 x_2 + q_2 x_3
	\end{pmatrix}
	.
\end{equation}
This is a Hamiltonian vector field, corresponding to the Hamiltonian $I_{x_1,x_2,x_3}\in C^{\infty}(T^{*}\R^{2})$ given by
\begin{equation}\label{eq:momentum_function}
	\begin{split}
		I_{(x_1,x_2,x_3)}(q_1,q_2,p_1,p_2) &= \frac{1}{2} \paren[\big]{q_1 q_2 + p_1 p_2 } x_1 \\
		&+ \frac{1}{2}\paren[\big]{ q_1 p_2  - q_2 p_1 }x_2 \\
		&+  \frac{1}{4}\Big(q_1^2 + p_1^2 - q_2^2 - p_2^2\Big) x_3.
	\end{split}
\end{equation}
By identifying $(\R^{3})^{*}$ with $\R^{3}$, it follows that the momentum map is given by $\pi$ in equation~\eqref{eq:extended_hopf}.


\section{Formulae for $X_{\pback{\hopf}{H}}$} 
\label{sec:lifted_vf}

In this appendix we give explicit coordinate formulae for $X_{\pback{\hopf}{H}}$ expressed in the derivatives of $H$.
Let $T^*_{\vect{z}}\hopf$ denote the \emph{cotangent lift}~\cite[\S\!~6.3]{MaRa1999} of $\hopf$ at the point $\vect{z}=(q_1,q_2,p_1,p_2)$, i.e., the transpose of the Jacobian matrix of $\hopf$:
\begin{equation}\label{eq:rectang_matrix}
	T^*_{\vect{z}}\hopf =
	\frac{1}{2}
	\begin{pmatrix}
		q_2 & \phantom{+}p_2 & \phantom{+}q_1 \\
		q_1 & -p_1 & -q_2 \\
		p_2 & -q_2 & \phantom{+}p_1 \\
		p_1 & \phantom{+}q_1 & -p_2
	\end{pmatrix}
	.
\end{equation}
Further, let $J$ be the canonical symplectic matrix on $T^*\R^{2}$.
Then the expression for $X_{\pback{\hopf}{H}}$ is
\begin{equation}\label{eq:expression_for_lifted_vf}
	X_{\pback{\hopf}{H}}(\vect{z}) = J^{-1}T^*_{\vect{z}}\hopf\,\ud{H}\big(\hopf(\vect{z})\big)
	= 
	\frac{1}{2}
\begin{pmatrix}
	\phantom{+}p_2 & -q_2 & \phantom{+}p_1 \\
	\phantom{+}p_1 & \phantom{+}q_1 & -p_2 \\
	-q_2 & -p_2 & -q_1 \\
	-q_1 & \phantom{+}p_1 & \phantom{+}q_2
\end{pmatrix}
\nabla H\big(\hopf(\vect{z}) \big)
	,
\end{equation}
which follows since $J X_{\pback{\hopf}{H}} = \ud \pback{\hopf}{H} = \pback{\hopf}{\ud H} = T^*\pi\circ\ud H\circ\pi$.

Notice that 
\begin{equation}\label{eq:collective_Hamiltonian_thm}
	X_{\pback{\hopf}{H}}(\vect{z}) = X_{I_{\ud H(\hopf(\vect{z}))}}(\vect{z}) = (\ud H(\hopf(\vect{z})))_{T^*\R^{2}}(\vect{z}),
\end{equation}
i.e., at each point $\vect{z}$ the vector field $X_{\pback{\hopf}{H}}$ is given by the infinitesimal generator corresponding to $\ud H(\hopf(\vect{z}))\in\su(2)$ evaluated at $\vect{z}$. (This result is called the \emph{Collective Hamiltonian Theorem} \cite[Th.~12.4.2]{MaRa1999}.)



\bibliographystyle{amsplainnat} 
\bibliography{collective}

\begin{thebibliography}{10}
\providecommand{\natexlab}[1]{#1}

\bibitem[{Aref(2007)}]{Ha2007}
H.~Aref, \emph{Point vortex dynamics: a classical mathematics playground}, J.
  Math. Phys. \textbf{48} (2007), 065401, 23.

\bibitem[{Guillemin and Sternberg(1984)}]{GuSt1984}
V.~Guillemin and S.~Sternberg, \emph{Symplectic techniques in physics},
  Cambridge University Press, Cambridge, 1984.

\bibitem[{Hairer et~al.(2006)Hairer, Lubich, and Wanner}]{HaLuWa2006}
E.~Hairer, C.~Lubich, and G.~Wanner, \emph{Geometric Numerical Integration},
  Springer-Verlag, Berlin, 2006.

\bibitem[{Libermann and Marle(1987)}]{LiMa1987}
P.~Libermann and C.-M. Marle, \emph{Symplectic geometry and analytical
  mechanics}, vol.~35 of \emph{Mathematics and its Applications}, D. Reidel
  Publishing Co., Dordrecht, translated from the French by Bertram Eugene
  Schwarzbach, 1987.

\bibitem[{Marsden and Weinstein(1983)}]{MaWe1983}
J.~Marsden and A.~Weinstein, \emph{Coadjoint orbits, vortices, and {C}lebsch
  variables for incompressible fluids}, Phys. D \textbf{7} (1983), 305--323,
  order in chaos (Los Alamos, N.M., 1982).

\bibitem[{Marsden and Ratiu(1999)}]{MaRa1999}
J.~E. Marsden and T.~S. Ratiu, \emph{Introduction to Mechanics and Symmetry},
  Springer-Verlag, New York, 1999.

\bibitem[{McLachlan et~al.(2013{\natexlab{a}})McLachlan, Modin, and
  Verdier}]{McMoVe2013b}
R.~I. McLachlan, K.~Modin, and O.~Verdier, \emph{{Collective Symplectic
  Integrators}},
  \href{http://arxiv.org/abs/1308.6620}{arxiv.org/abs/1308.6620},
  2013{\natexlab{a}}.

\bibitem[{McLachlan et~al.(2013{\natexlab{b}})McLachlan, Modin, Verdier, and
  Wilkins}]{McMoVeWi2013}
R.~I. McLachlan, K.~Modin, O.~Verdier, and M.~Wilkins, \emph{Geometric
  generalisations of {S}hake and {R}attle}, Foundations of Computational
  Mathematics  (2013{\natexlab{b}}).

\bibitem[{McLachlan and Quispel(2002)}]{McQu2002}
R.~I. McLachlan and G.~R.~W. Quispel, \emph{Splitting methods}, Acta Numer.
  \textbf{11} (2002), 341--434.

\bibitem[{Weinstein(1983)}]{We1983}
A.~Weinstein, \emph{The local structure of {P}oisson manifolds},
  J.~Differential Geom. \textbf{18} (1983), 523--557.

\end{thebibliography}


\end{document}